\documentclass[12pt,twoside,a4paper,final]{article}
\usepackage{latexsym,amssymb}
\title{Central limit theorems for Coxeter systems\\
and Artin systems of extra large type}
\author{Gero Fendler \\ 
Naturwissenschaftlich Technische Fakult\"{a}t I\\
Fachrichtung 6.1 Mathematik \\
Universit\"at des Saarlandes \\ 
}
\newcounter{prop}
\newcounter{def}

\newtheorem{thm}{Theorem}

\newtheorem{lem}{Lemma}

\newtheorem{prop}[prop]{Proposition}

\newenvironment{proof}[1][Proof:]{\it #1 \rm }{\qed}
\newenvironment{proofth}[1]{{\it Proof of Theorem #1:}\\}{\qed}
\newcommand{\NN}{\ensuremath{\mathbb{N}}}
\newcommand{\CC}{\ensuremath{\mathbb{C}}}
\newcommand{\id}{\ensuremath{\mbox{\rm Id}}}
\newcommand{\word}[3]{(#1_{#2} ,\ldots ,#1_{#3})}
\newcommand{\E}{\mbox{E}}
\newcommand{\var}{\mbox{V}}
\newcommand{\la}[1]{{\rm \bf l}(#1)}
\newcommand{\eval}[1][.]{\mathbf{ev}(#1)}
\newcommand{\card}[1][\#]{\protect{#1}}
\newcommand{\qed}{\nopagebreak \hfill \ensuremath{ \underline{~~~~\Box}}\\}
\newcommand{\NC}{{\rm NC}}
\newcommand{\prob}[1]{\ensuremath{\mathbf{prob}}(#1)}
\newcommand{\abs}[1][.]{|{#1}|}
\newlength{\arrowlength}%
\newcommand{\notconv}{%
\settowidth{\arrowlength}{ \ensuremath{\displaystyle{\longrightarrow}}}%
\makebox[0pt][l]{\hspace{.10\arrowlength}
\ensuremath{\displaystyle{\backslash}}}
\ensuremath{\displaystyle{\longrightarrow}}}
\begin{document}
\maketitle
\thanks{{\bf Acknowledgements :} The author is grateful to M. Bo\.{z}ejko and
R. Speicher for useful discussions on the theme of this note. The idea 
to consider randomly chosen Coxeter groups came up in a discussion with
R. Speicher.}
\par
\begin{abstract}
Given a Coxeter system of large type we prove a non--commutative 
central limit theorem:\\
After normalisation with the square root of $n$ the characteristic 
function of the set of the first $n$ generators tends in distribution 
to Wigners semi--circle law.
\par
If one chooses the group presentation at random, now admitting a certain
amount of commutation among the generators, the distributional limit 
will be the $q$--deformed semi--circle law, $1\leq q \leq 1$, almost 
surely.
\end{abstract}
\section{Introduction}%
\label{sec:intro}
We should start by explaining our understanding of an infinite
Coxeter system.
Loosely speaking this is a presentation of an infinitely 
generated group of a specific type.
\\
Let $S$ be a countably infinite set and for $s,s' \in S$ let
$m(s,s') \in \NN \bigcup \{\infty\}$ be given. 
Assume $m(s,s)=1 $ and $m(s,s')\geq2$ for $s\neq s'$. 
Then the group $G_m$ corresponding to these
data has the presentation:
\begin{eqnarray}
\label{presentation1}
(ss')^{m(s,s')} & = & \id, \quad  \mbox{ for }s,s' \in S \mbox{ with }
m(s,s') \ < \ \infty.
\end{eqnarray}
If $J \subset J'$ are two subsets of $S$ then, taking $m_J$ and $m_{J'}$
as the sub-matrices of $m$ obtained by restricting to $J$ resp. $J'$ as
sets of generators, Coxeter groups
$W_J$ and $W_{J'}$ are defined by the corresponding presentations. Moreover
the inclusion $J \subset J'$ identifies $W_J$ to a parabolic subgroup
of $W_{J'}$. We refer to the nice exposition of Humphreys~\cite{Hum_90}
for further information on parabolic subgroups.
As the above group $G_m$ we take the inductive limit,
along the net of finite subsets of $S$, with these identifications. 
Then it is clear that the above
relations (\ref{presentation1}) 
are fulfilled in $G_m$. Further, since a relation in $G_m$ involves
only finitely many generators, it really  is a relation in some finitely generated
subgroup $W_J$. This shows that any relation in $G_m$ is a consequence of  
(\ref{presentation1}).
\par
Now we consider $\CC G_m$, in the natural way,
\begin{eqnarray}
\label{eq:conv}
\sum_{x\in G_m}a_x \delta (x) \times \sum_{x\in G_m}b_x \delta (x)%
&\mapsto& \sum_{z\in G_m} \sum_{\{x,y \in G_m : xy=z\}} a_x b_y \delta (z)
\end{eqnarray}
as a convolution algebra and denote 
\begin{eqnarray}
\label{eq:trace1}
\varphi : \CC G_m& \rightarrow & \CC \\
\label{eq:trace2}
\sum_{x\in G_m}a_x \delta (x) & \mapsto & a_{\id}
\end{eqnarray}
the canonical trace on this algebra.
\par
An algebraic {\em central limit theorem} in this setting reads:
\begin{thm}
\label{coxclt}
Let $(G_m,S)$ be a Coxeter system with a countable arbitrarily enumerated  set
$S=\{s_1, s_2, s_3, \ldots \}$ of generators.
Assume, that 
$m(s,s') \geq 3 \mbox{ for all } s\not=s'\;s,s' \in S $.
Then for all $k \in \NN$
\[
\lim_{N\rightarrow \infty}\ %
 \varphi \left( \left(\frac{\delta(s_1)+ \ldots + \delta(s_N)}{\sqrt{N}}%
\right)^k\right)%
= \left\{
\begin{array}{cl}
0 &\mbox{ if } k \mbox{ is odd}\\
\displaystyle{\frac{1}{n+1}} \left(\!\!
\begin{array}{c}
2n \\
n
\end{array}\!\!\right) & 
\mbox{ if } k=2n \mbox{ is even}.
\end{array}
\right.
\]
\end{thm}
Since the Catalan numbers $\frac{1}{n+1}\left(\!\!\begin{array}{c}
2n \\
n
\end{array}\!\!\right)$ are just the $2n$--th moments of Wigners semicircular law
$d\mu(t)=\frac{1}{2\pi}\chi_{[-2,2]}(t)\sqrt{4-t^2}\,dt$ we
recover it, at least in the sense of convergence of moments, as a central limit.
\par
To consider central limit theorems of this type is inspired 
from papers of Bo\.{z}ejko and Speicher, see \cite{bozspei91},
\cite{bozspei92}. The probabilistic 
approach to $q$-interpolated limit theorems follows the methods
of Speicher \cite{speicher92}.
\par
In section~\ref{sec:artin} we consider the corresponding problem
for Artin groups of extra-large type, using results of Appel and Schupp
\cite{appelschupp83}.

\section{Proofs for Coxeter groups}
\label{sec:cox}
We shall be interested in the manner in which words in the language generated
by $S$ reduce or do not reduce in $G_m$. First we recall.
\begin{lem}
\label{lem:1}
Assume $I \subset S$ and $s\notin I$. If $w\in W_I$ has the reduced expansion
$w=w_1 \ldots w_n$, then $w_i \in I, \mbox{ for } 1 \leq i \leq n$,
 and $ws=w_1 \ldots w_n s$ and
$sw=sw_1 \ldots w_n$ are both reduced.
\end{lem}
\begin{proof}
The first assertion is contained in part (b) of the theorem in 
section 5.5 of~\cite{Hum_90}. The others then follow immediately 
from the exchange condition, c.f.\ sec.~5.8~of \cite{Hum_90}.
\par
For example, if $w=w_1 \ldots w_n s$ would not be reduced, then
for a unique $k\in \{1,\ldots, n \}$ 
\[ws=w_1 \ldots \hat{w_k} \ldots w_n \]
(omitting $\hat{w_k}$). Hence
\[
w_{k}\ldots w_n s = w_{k+1} \ldots w_n.
\]
Or
\[
s=w_n \ldots w_{k+1} w_k w_{k+1} \ldots w_n.
\]
By the deletion condition, see the Corollary in 5.8 of \cite{Hum_90},
a reduced expansion of the right hand side may be obtained by 
deleting pairs of letters. But this implies $s=w_j$ for some $j$, 
in contradiction to $s\notin I$.
\par
The above representation of $sw$ is dealt with analogously.
\end{proof}
\begin{lem}
\label{lem:2}
Let $\word{w}{1}{r}$ be a word in the generators.
Assume that $s\in S$ appears only once among the above letters. Then
\begin{description}
\item[~(i)]
\[ w \ = \ w_1\cdot \ldots \cdot w_r \ \neq \id \quad \mbox{in} \quad G_m.\]
\item[~(ii)] In each reduced representation of this group element
\[w \ = \ t_1\cdot \ldots \cdot t_l , \ t_i \in S \]
there appears $s$, i.e.\ there exists a 
$j \in \{1,\ldots,l\}$  
with $s=t_j$.
\end{description}
\end{lem}
\begin{proof}
Assume $s=w_i$, 
$s \notin I = \{ w_1, \ldots , w_{i-1},w_{i+1}, \ldots , w_r\}$.
Then $u=w_1\cdot \ldots \cdot w_{i-1}$ and 
$v =w_{i+1}\cdot \ldots \cdot w_r$
are elements of $W_I$, and $\id = w =usv$ would imply 
$s = u^{-1}v^{-1} \in W_I$.
Hence $s \in I $, a contradiction.
\par
If we would have $w=t_1\cdot \ldots \cdot t_l$, with 
$s\notin J =\{t_1,\ldots , t_l\}$
then $w=usv \in W_J$ hence, as above, $s \in W_{(I \cup J)}$,
contradicting $s \notin I \cup J$. 
\end{proof}
\par
To a word $\word{w}{1}{r}$ there is associated a partition
$V= \{ V_1, \ldots, V_p\}$ of the index set $\{1, \ldots , r\}$ by
taking just the pre-images of single points under the map
$i \mapsto w_i$ as a map defined on $\{1, \ldots, r\}$ to
$S$.
\par
We recall
\begin{description}
\item[~~(i)]
A partition $V= \{ V_1, \ldots, V_p\}$ is called a pair-partition
if all its elements are two element sets.
\item[~(ii)]
A partition $V= \{ V_1, \ldots, V_p\}$ of $\{1, \ldots , r\}$  is 
called crossing if there are $n\neq m$ and
$i<k<j$, and $i,j \in V_m$, $ k\in V_n$.
The partition is called non-crossing otherwise.
\item[(iii)]
For a pair-partition the condition of being non-crossing amounts to:
$i<k<j$, and $i,j \in V_m$, $ k\in V_n$ imply that 
$V_n \subset \{i+1, \ldots , j-1\}$.
\end{description}
\begin{lem}
\label{lem:3}
Assume that for all $s,s' \in S$ we have $m(s,s') \geq 3$.
If $\word{w}{1}{r}$ is a word which in the above sense defines a
pair-partition $V$, then
\begin{eqnarray}
w\ = \ w_1\cdot \ldots \cdot w_r \ = \ \id & \mbox{if and only if} &
V \mbox{ is non-crossing}.
\end{eqnarray}
\end{lem}
\begin{proof}
Assume that a word $\word{w}{1}{r}$
defines a crossing partition and that
$w\ = \ w_1\cdot \ldots \cdot w_r \ = \ \id $.
Let $s=w_i=w_j \neq t=w_k=w_l$ with
$i<k<j<l$ and denote $I=\{w_n:n\notin\{i,j,k,l\}\}$.
We write $w = usxtysztv$, where $u,v,x,y,z \in W_I$ and let
$xty=t_1\ldots t_n$ be a reduced expansion.
From Lemma~\ref{lem:1} we know that $sxty=st_1\ldots t_n$ and
$xtys=t_1\ldots t_ns$ are reduced expansions.
\par
Moreover $sxtys=st_1\ldots t_ns$ is a reduced expansion too.
For if we had $\la{sxtys}<\la{sxty}$ then by Corollary~1 of \cite{szwarc_98}
it would follow that in the above reduced expansion of $sxty$ the last or
last but one letter equals $s$. Since $t$ must appear among 
$t_1, \ldots, t_n$ we would have $sxty=st$ and hence 
$sts=u^{-1}v^{-1}tz^{-1} \in W_{I\cup\{t\}}$. 
Here the left hand side could not be
reduced, since $s\notin I\cup\{t\}$. But then 
by the deletion condition,~c.f.\ sec.5.8 of \cite{Hum_90}, 
only $sts=t$ would be possible, i.e.\ $m(s,t)=2$.
\par
Now, again from $w=\id$, it follows that 
$st_1\ldots t_ns=sxtys= u^{-1}v^{-1}tz^{-1}\in W_{I\cup\{t\}}$.
In contradiction to $s\notin I\cup\{t\}$ and the fact that left hand side
is reduced.
\par
The converse implication is proved by a straightforward induction.
\end{proof}
\mbox{}\\
\begin{proofth}{\ref{coxclt}} 
Denote 
$\mathcal{V}_p^k=\{V : V \mbox{ is a partition of} \{1, \ldots , k\}  
\mbox{ into }p\mbox{ sets } \}$
and for a partition $V\in \mathcal{V}_p^k $ let
$\mathcal{W}_{V}^N=\{(s_{j_1}, \ldots , s_{j_k}) : \mbox{  the word }%
\word{s}{j_1}{j_k}%
\mbox{ defines } V , 1\leq j_i \leq N \mbox{ for } i=1,\ldots, k\}$.
We compute:
\begin{eqnarray}
\label{eqn:c1}
\varphi \left( \left(\frac{\delta(s_1)+ ... + \delta(s_N)}{\sqrt{N}}%
\right)^k\right)& = %
&\sum_{i_1, \ldots, i_k =1}^N \left(\frac{1}{\sqrt{N}}\right)^k%
\varphi \left(s_{i_1}\cdot \ldots \cdot s_{i_k}\right)\\
\label{eqn:c2}
&=& \sum_{p=1}^k N^{-\frac{k}{2}}%
\sum_{V\in \mathcal{V}_p^k}%
\sum_{w\in \mathcal{W}_{V}^N} \varphi ( \eval[w])
\end{eqnarray}
Here $\eval$ sending $w = \word{s}{i_1}{i_m}$ to the product
$s_{i_1}\cdot \ldots \cdot s_{i_m}$
denotes the evaluation map from the set of words in the
generators to group elements.
\par 
Now, by the deletion condition $\eval[\word{s}{i_1}{i_k}] \neq \id$,
whenever $k$ is odd. The assertion of the theorem being established then.
We may henceforth assume that $k=2n$ is even.
\par
If a partition $V$ of the $2n$ element set $\{1,\ldots,k\}$
contains more than $n$ of its (pairwise disjoint) subsets
then it must contain a one element set and therefore,
by Lemma~\ref{lem:2}, $\varphi(\eval[w])=0$  whenever
$w\in \mathcal{W}_{V}^N$ for some $V\in \mathcal{V}_p^k$ with $p>n$.
The above sum reducing to:
\begin{eqnarray}
\label{eqn:c3}&&
\sum_{p=1}^n N^{-\frac{k}{2}}%
\sum_{V\in \mathcal{V}_p^k}%
\sum_{w\in \mathcal{W}_{V}^N} \varphi ( \eval[w]).
\end{eqnarray}
Taking into account
that there are $A_{N,p} = N (N-1) \cdot \ldots \cdot (N-p+1)$ 
words in the letters
$s_1, \ldots , s_N$ which define a partition $V\in \mathcal{V}^k_p$,
this just equals 
\begin{eqnarray}
\label{eqn:c4}&&%
N^{-\frac{k}{2}}%
\sum_{V\in \mathcal{V}_n^k}%
\sum_{w\in \mathcal{W}_{V}^N} \varphi ( \eval[w]) + o(1).
\end{eqnarray}
Any partition of a $2n$ element set in $n$ subsets not containing a
one element set must be a pair-partition. 
By Lemma~\ref{lem:3} those which contribute to the sum are exactly the 
non-crossing ones. Hence we end with
\[
\card{\NC_2(k)}%
\frac{N}{N}\frac{N-1}{N}\cdot \ldots \cdot\frac{N-\frac{k}{2}+1}{N}%
 + o(1),
\]
where $\NC_2(k)$ denotes the set of non-crossing pair-partitions
of a k element set.
This cardinality has been computed in~\cite{Kreweras72} to be
the Catalan numbers which finishes our proof.
\end{proofth}
\section{Probabilistic Interpolation}%
\label{sec:prob}
As we have seen in Theorem~\ref{coxclt} we obtain in the
limit the moments of the semicircle law, whenever
there is no commutation at all in the generators of the Coxeter system.
On the other hand it is not difficult to compute
the limit measure, when all generators commute.
In fact we are then treating independent Bernoulli random variables 
and proving the classical DeMoivre-Laplace theorem.
\par
The aim of this section is to interpolate between these situations
by randomly choosing the Coxeter system. Thus let
$S=\{s_1,s_2,\ldots\}$ be an enumerated infinite generating set.
We shall consider a Coxeter matrix $m=\left(m(s,t)\right)_{s,t\in S}$
as a random variable satisfying the following
independence conditions and requirements on the distribution:
\begin{description}
\item[(i)]%
Of course $m(s,t)=m(t,s)$ and $m(s,s)=1$.
\item[(ii)]%
If $s\neq t$ and  $s'\neq t'$, and neither $(s,t)= (s',t')$
nor $(s,t)=(t',s')$, as ordered pairs, then 
$m(s,t)$ and $m(s',t')$ are independent, identically
distributed.
\item[(iii)]%
For some $p\in [0,1]$ for all $s\neq t$:
\[
\prob{m(s,t)=2}=p \mbox{ and } \prob{ m(s,t)\in \{3,4,\ldots,\infty \} }=1-p.
\] 
\end{description}
For our topic we now need
a refinement of Lemma~\ref{lem:3}.
\begin{lem}
\label{lem:p1}%
Let $G_m,S$ be a Coxeter system with Coxeter matrix
$\left( m(s,s') \right)_{s,s' \in S}$.
Assume that the word $\word{w}{1}{r}$ defines a pair-partition.
If there are $i<k<j<l$ and $s,t\in S$ such that
$s=w_i=w_j$ and $t=w_k=w_l$, 
so that the partition contains at least this specified crossing,
then 
\begin{eqnarray}
w \ = \ w_1\cdot \ldots \cdot w_r \ = \ \id &\mbox{implies}& m(s,t)\ = \ 2.
\end{eqnarray}
\end{lem}
\begin{proof}
As in the proof of Lemma~\ref{lem:3} we let $I=\{w_n:n\notin\{i,j,k,l\}\}$
and write $w=usxtysztv$, with $u,x,y,z,v \in W_I$.
Let again $xty=t_1 \ldots t_n$ be a reduced expansion.
\par
We claim that $st_1 \ldots t_ns$ is reduced, except if
$m(t_r,s)=2$ for $r=1,\ldots,n$.
Now, if $st_1 \ldots t_ns$ is not reduced then by the deletion condition
either
\[
\quad t_1 \ldots t_n = st_1 \ldots t_ns
\]
or $\mbox{ for some } m'\in \{1, \ldots , n\}$:
\[
\quad st_1 \ldots \hat{t}_{m'}\ldots t_n \hat{s} = st_1 \ldots t_ns%
\]
The second would imply
$t_{m'}t_{m'+1}\ldots t_n s $ and a forteriori $t_1\ldots t_n s$
not being reduced, in contradiction to
 Lemma~\ref{lem:1}.
\par
We are left with
\[
t_1 \ldots t_n s = st_1 \ldots t_n
\]
Here both sides are reduced by Lemma~\ref{lem:1}.
A Lemma of Deodhar, see Proposition 1 of \cite{szwarc_98}
(taking the notation from that paper), implies that there is 
a reduced expansion:
\[
t_1 \ldots t_n = N(a_m,b_m)\ldots N(a_1,b_1)
\]
where $c(a_r,b_r)=a_{r+1}, d(a_r,b_r) \neq b_{r+1}$ and 
$a_1=s$, $t_n=b_1$, $c(a_m,b_m)=s$.
Here $s$ would appear as a letter on the right hand side except
if $m(s,b_1)=m(s,t_n)=2$ and $a_m=s$,  $m(s,b_m)=2$. 
But then $s=c(a_{m-1},b_{m-1})$, and again $s$ would appear except
if $m(a_{m-1},b_{m-1})=2$, $a_{m-1}=s$. Inductively we obtain
$a_m=a_{m-1}= \ldots =a_1=s$ and $m(s,b_r)=2$ for $r=1,\ldots ,m$.
Since $\{t_1,\ldots,t_n\}=\{b_m,\ldots,b_1\}$ we have established the claim.
\par
Finally, $t$ must appear among $t_1,\ldots,t_n$ and the proof can be completed
exactly as in the proof of Lemma~\ref{lem:3}.
\end{proof}
\par
Given a pair-partition $V=\{V_1,\ldots,V_r\}$
of the set $\{1, \ldots, 2r\}$ the sets $V_i=\{e_i,f_i\}$
are in fact naturally ordered, we shall assume $e_i<f_i$.
By renaming the subsets of the Partition we may assume
further, that $e_i<e_k$ if $i<k$.
The set of inversions  
of the partition is then defined
as
\[I(V)=\{(i,j):e_i<e_j<f_i<f_j\}.\]
An example is visualised in figure~1.
\begin{figure}
\label{fig:1}
\caption{}
\setlength{\unitlength}{3947sp}%
\begingroup\makeatletter\ifx\SetFigFont\undefined%
\gdef\SetFigFont#1#2#3#4#5{%
  \reset@font\fontsize{#1}{#2pt}%
  \fontfamily{#3}\fontseries{#4}\fontshape{#5}%
  \selectfont}%
\fi\endgroup%
\begin{picture}(5160,1460)(346,-634)
\thicklines
\put(1126,-420){\oval(1350,1350)[tr]}
\put(1126,-420){\oval(1350,1350)[tl]}
\put(1126,-420){\oval(450,450)[tr]}
\put(1126,-420){\oval(450,450)[tl]}
\put(3830,-420){\oval(3142,2240)[tr]}
\put(3830,-420){\oval(3142,2240)[tl]}
\put(3376,-420){\oval(1348,1288)[tr]}
\put(3376,-420){\oval(1348,1288)[tl]}
\put(4031,-420){\oval(1800,1800)[tr]}
\put(4031,-420){\oval(1800,1800)[tl]}
\put(4051,-420){\oval(900,900)[tr]}
\put(4051,-420){\oval(900,900)[tl]}
\put(451,-530){\circle{180}}
\put(882,-530){\circle{180}}
\put(1351,-530){\circle{180}}
\put(1806,-530){\circle{180}}
\put(2251,-530){\circle{180}}
\put(2701,-530){\circle{180}}
\put(3132,-530){\circle{180}}
\put(3606,-530){\circle{180}}
\put(4051,-530){\circle{180}}
\put(4501,-530){\circle{180}}
\put(4951,-530){\circle{180}}
\put(5401,-530){\circle{180}}
\end{picture}
\end{figure}
\par
When a word $w=\word{w}{1}{2r}$ defines a pair-partition $V$,
then each of the sets $V_i$ of $V$
bears a label $t_i(w) \in S$, namely the image of $V_i$
under the map $V_i\mapsto w_{e_i}=w_{f_i} \in S$.
To an element $(i,j) \in I(V)$ we may associate
this way a pair $(t_i(w),t_j(w))$ of elements of $S$
and hence the number $m_{i,j}(w) := m(t_i(w),t_j(w))$,
where $m$ is the Coxeter matrix.
Reformulating Lemma~\ref{lem:p1} and taking Lemma~\ref{lem:3} into account:
\begin{lem}
\label{lem:p2}
Assume, that the  word $w=\word{w}{1}{2r}$ defines a pair-partition $V$.
Then
\begin{eqnarray*}
\eval[w]&=&\id
\end{eqnarray*}
if and only if for all $(i,j) \in I(V)$: 
\begin{eqnarray*}
m_{i,j}(w)&=&2.
\end{eqnarray*}
\end{lem}
\begin{prop}
\label{prop:p1}
For a pair-partition $V$ of $\{1,\ldots,2r\}$ consider the random variable
\[
X_N = \frac{1}{N^r} \sum_{w\in \mathcal{W}_{V}^N} \varphi(\eval[w]). 
\]
Then, almost surely
\[
X_N \rightarrow p^{\card{I(V))}}
\]
as $N \rightarrow \infty$.
\end{prop}
\begin{proof}
We first compute the expectation of $X_N$.
Let $w=\word{w}{1}{2r}$ be a word in the letters $s_1, \ldots ,s_N$.
If $(i,j)$ and $(k,l)$
are different inversions of $V$, then $m_{i,j}(w)$ and $m_{k,l}(w)$
are independent random variables. Hence, by {\bf (ii)} and {\bf (iii)},
\[
\prob{m_{i,j}(w)=2 , \mbox{ for all } (i,j)\in I(V)} = p^{\card{I(V)}}.
\]
Since by Lemma~\ref{lem:p2} $\varphi(\eval[w])=1$
exactly if $m_{i,j}(w)=2 , \mbox{ for all } (i,j)\in I(V)$ we obtain
\begin{eqnarray*}
\E(X_N)&=&\frac{1}{N^r} \sum_{w\in \mathcal{W}_{V}^N} \E(\varphi(\eval[w]))\\
&=&\frac{1}{N^r} \sum_{w\in \mathcal{W}_{V}^N} p^{\card{I(V)}}\\
&=&\frac{N(N-1)\ldots(N-r+1)}{N^r} p^{\card{I(V)}}
\end{eqnarray*}
which converges to $p^{\card{I(V)}}$.
\par
Next we compute the variance:
\begin{eqnarray*} 
\var(X_N)&=&\E(X_{N}^2)-E(X_N)^2\\
&=&\frac{1}{N^{2r}}\left\{\sum_{w\in \mathcal{W}_{V}^N}\sum_{v\in \mathcal{W}_{V}^N}%
\E(\varphi(\eval[w])\varphi(\eval[v])) -%
\left(\frac{N!}{(N-r)!}\right)^2 p^{2\card{I(V)}}\right\}\\
&=&\frac{1}{N^{2r}}\sum_{w\in \mathcal{W}_{V}^N}\sum_{v\in \mathcal{W}_{V}^N}%
\left\{ \E(\varphi(\eval[w])\varphi(\eval[v])) -%
p^{2\card{I(V)}}\right\}.
\end{eqnarray*}
If $\varphi(\eval[w])$ and $\varphi(\eval[v])$ are independent, then the
corresponding summand vanishes.
But, if these random variables are dependent, then
for some $(i,j) \in I(V)$ we have 
\begin{eqnarray*}
(t_i(w),t_j(w))&=& (t_i(v),t_j(v)).
\end{eqnarray*}
For given 
$w\in \mathcal{W}_{V}^N$ and $(i,j) \in I(V)$ there are at most 
$(N-2)(N-3)\ldots (N-r+1)$ elements 
$v\in \mathcal{W}_{V}^N$ fulfilling this equality.
Since  $\card{\mathcal{W}_{V}^N}=N(N-1)\ldots (N-r+1)$, we
conclude that for at most 
\[\card{I(V)}\ {N(N-1)\left((N-2)(N-3)\ldots (N-r+1)\right)^2}\]
pairs of elements the random variables in question are dependent.
In this case the contribution of a summand is at most $1-p^{2\card{I(V)}}$.
hence we may estimate:
\begin{eqnarray*}
\var(X_N)&\leq&\card{I(V)}\frac{1}{N^2}.
\end{eqnarray*}
As is well known,
the summability of $\sum_{N=1}^{\infty}\var(X_N)$
now implies the almost sure convergence:\\
In the underlying probability space $(\Omega,\prob)$ let 
\[A_{n,k} =\bigcup_{N>n} \{\abs[X_N-\E(X_N)]>\frac{1}{k}\}.\]
Then the set, where $X_N$ does not converge to 
$x^0=\lim_{N\rightarrow \infty} \E(X_N)$
can be written as
\begin{eqnarray*}
\{X_N \notconv x^0\}&=& \{\abs[X_N-\E(X_N)] \notconv 0 \}\\
&=&\bigcup_{k\in \NN} \bigcap_{n\in \NN}A_{n,k}.
\end{eqnarray*}  
Since the sets $A_{n,k}$ decrease in $n$ it is sufficient to
show for each $k$ that their probability tends to $0$ as n
tends to infinity. \\
Clearly
$A_{n,k} = \{\sup_{N>n}\abs[X_N-\E(X_N)]>\frac{1}{k}\}$.
Hence, by Chebycheff's inequality:
\begin{eqnarray*}
\prob{A_{n,k}}&=&\prob{\{\sup_{N>n}\abs[X_N-\E(X_N)]>\frac{1}{k}\}}\\
&\leq&{k^2}\E(\sup_{N>n}\abs[X_N-\E(X_N)]^2)\\
&\leq&{k^2}\left(\sum_{N=n}^{\infty}\abs[X_N-\E(X_N)]^2\right).
\end{eqnarray*}
\end{proof}
We conclude this section with its main theorem:
\begin{thm}
\label{probCLT}
Let $(G_m,S)$ be a random Coxeter system with a countable arbitrarily enumerated  set
$S=\{s_1, s_2, s_3, \ldots \}$ of generators, where we assume that the 
random Coxeter matrix $m$ fulfils the conditions {\bf (i), (ii)}
and {\bf (iii)}.

Then, almost surely,  for all $k \in \NN$
\[
\lim_{N\rightarrow \infty}\ %
 \varphi \left( \left(\frac{\delta(s_1)+ \ldots + \delta(s_N)}{\sqrt{N}}%
\right)^k\right)%
= \left\{
\begin{array}{cl}
0 &\mbox{ if } k \mbox{ is odd}\\
\displaystyle{\sum_{V\in \mathcal{V}_r^k,\,V\mbox{ a pair partition}}%
p^{\card{I(V)}}}& 
\mbox{ if } k=2r \mbox{ is even}.
\end{array}
\right.
\]
\end{thm}
\begin{proof}
As in the prove of Theorem~\ref{coxclt} we see that for even $k=2r$:
\begin{eqnarray*}
\sum_{i_1, \ldots, i_k =1}^N \left(\frac{1}{\sqrt{N}}\right)^k%
\varphi \left(s_{i_1}\cdot \ldots \cdot s_{i_k}\right)&=&%
N^{-r}%
\sum_{V\in \mathcal{V}_r^k}%
\sum_{w\in \mathcal{W}_{V}^N} \varphi ( \eval[w]) + o(1),
\end{eqnarray*}
whereas the left hand side vanishes
for odd $k$.
\par
By our Proposition~\ref{prop:p1} this, almost surely tends to
\begin{eqnarray*}
\lim_{N\rightarrow \infty}%
\sum_{V\in \mathcal{V}_r^k} X_N &=& \sum_{V\in \mathcal{V}_r^k}%
p^{\card{I(V)}}.
\end{eqnarray*}
\end{proof}

\section{Artin Groups}%
\label{sec:artin}
Let $S$ be  a finite set and $m$ a Coxeter matrix over $S$.
The Artin group $A$ corresponding to $m$ then is the group
with generating set $\{a_s:s \in S \}$ and defining relations, given for
$s \neq s' , \ s,s' \in S$ with $m(s,s')<\infty$, by
\[ a_sa_{s'}a_s \ldots = a_{s'}a_s a_{s'} \ldots \ , \]
where both products have $m(s,s')$ factors.
\par
The corresponding Coxeter group $G_m$ had the additional relations
\[ a_s^2= \id .\]
The map $a_s \mapsto s,\ s\in S$ from the Artin group $A$ to the Coxeter
group $G_m$ hence extends to a homomorphism $\Phi$ containing the
normal subgroup $N$ generated by the set $\{ a_s^2: s \in S\}$ in its kernel.
On the other hand all Coxeter relations are fulfilled by the cosets
$\mbox{mod } N$ of the generators of $A$.
Thus we may identify 
the quotient
$A/N$ with $G_m$.
\par
In that case of Coxeter groups we had for $I \subset S$, 
still $S$ a finite set, an
isomorphism of the Coxeter group obtained from the restriction of
$m$ to $I\times I$ and the subgroup $W_I$ of $G_m$ generated by the set
$I$. This enabled us to define, for an infinite set $S$ the corresponding
Coxeter group as an inductive limit along the net of finite subsets of $S$.
For general Artin groups I don't know whether this is possible. But for
the case that $m(s,s')\geq 3$ whenever $s\neq s'$ Appel and Schupp 
\cite{appelschupp83} showed that this is indeed true, cf. their
Corollary~3.
Furthermore, if $m(s,s')\geq 4$ whenever $s\neq s'$,
they call those
Artin groups of extra large type,
then a word, reduced in $A$,   represents an element of a subgroup
$A_J$ generated by the set $\{a_j: j \in J\}$ only if it is a
word on this set.
 and we shall in this section discuss 
Artin groups of extra large type with countable infinite generating set
$S$ as the inductive limit along the finite subsets
defined as for Coxeter groups.   
\par
As in (\ref{eq:conv}) we consider $\CC A$ as a convolution algebra and
$\varphi_A: \CC A \rightarrow \CC$ defined as in (\ref{eq:trace2})
the canonical trace. Let $S=\{s_1, s_2,s_3, \ldots \}$ be enumerated
and denote $a_i=a_{s_i}$ $i \in \NN$ the generators of the 
Artin group.
\begin{thm}
\label{artinclt}
Under the above conditions we have
\[
\lim_{N\rightarrow \infty}\ %
 \varphi_A \left( \left(\frac{\delta(a_1)+\delta(a_1^{-1})+ \ldots + %
\delta(a_N)+\delta(a_N^{-1})}{\sqrt{2N}}%
\right)^k\right)%
= \left\{
\begin{array}{cl}
0 &\mbox{if } k \mbox{ is odd}\\
\displaystyle{\frac{1}{n+1}} \left(\!\!
\begin{array}{c}
2n \\
n
\end{array}\!\!\right) & 
\mbox{if } k=2n \mbox{ is even}.
\end{array}
\right.
\]
\end{thm}
\par
Given a word $(a_{i_1}^{\epsilon_1}, \ldots , a_{i_k}^{\epsilon_k})$, 
$\epsilon_i \in \{+1,-1\}$ for $i=1,\ldots , k$,
in the generators and their inverses we associate a partition 
$V$ of $\{1,\ldots,k\}$ this time by $V=\{ V_1, \ldots ,V_p\}$,
where $p=\card\{i_1,\ldots, i_k\}$, by taking the pre-images of the 
sets $\{j\}, \ j \in \NN$, under the map $ k \mapsto i_k$.
Clearly the evaluation in the Artin group can yield the identity
at most if $(\Phi{a_{i_1}^{\epsilon_1}},\ldots,\Phi{a_{i_1}^{\epsilon_k}})$
evaluates to the identity element in the Coxeter group. Hence we
obtain from Lemma~\ref{lem:2} and \ref{lem:3}:
\begin{lem}
\label{lem:4}
Let $w=(a_{i_1}^{\epsilon_1}, \ldots , a_{i_k}^{\epsilon_k})$ be a 
word in the generators and its inverses with associated partition $V$.
\begin{description}
\item[~(i)]
If $V$ contains a one-element set, then
\[ \eval[a_{i_1}^{\epsilon_1}, \ldots , a_{i_k}^{\epsilon_k}] \neq \id.\]
\item[(ii)]
If $V$ is a pair-partition, then 
\[ \eval[a_{i_1}^{\epsilon_1}, \ldots , a_{i_k}^{\epsilon_k}] = \id.\]
implies that $V$ is non-crossing.
\end{description}
\end{lem}
\par
To prove the theorem we need one more Lemma:
\begin{lem}
\label{lem:5}
Assume that $w=(a_{i_1}^{\epsilon_1}, \ldots , a_{i_k}^{\epsilon_k})$
defines a pair-partition. 
If for some 
$j,l\in \{1,\ldots,k\}$ with $j\neq l$ we have $a_{i_j} = a_{i_l}$ and
$\epsilon_j=\epsilon_l$ then
\[ \eval[a_{i_1}^{\epsilon_1}, \ldots , a_{i_k}^{\epsilon_k}] \neq \id.\]
Conversely, if $w$ defines a non-crossing pair-partition and
$\epsilon_j=-\epsilon_l$ whenever $a_{i_j} = a_{i_l}$, $j\neq l, \ j,l\in
\{1,\ldots, k\}$ then $\eval[w]=\id$.
\end{lem}
\begin{proof}
We  assume 
$\eval[a_{i_1}^{\epsilon_1}, \ldots , a_{i_k}^{\epsilon_k}] = \id$
and derive a contradiction.
Possibly taking the inverse we can suppose $\epsilon_j=\epsilon_l=1$,
and of course $j<l$ too. Let $I= \{ i_{j+1}, \ldots, i_{l-1} \}$ and
$J=\{i_{r} : r<j \mbox{ or } l<r\}$.
Then, for some $s \notin I\cup J$,
\[\eval[a_{i_1}^{\epsilon_1}, \ldots , a_{i_k}^{\epsilon_k}] = ua_s w a_s v\]
with $u,v \in A_J$ and $w\in A_I$. Here for $K\subset S$ 
$A_K$ is the subgroup of $A$
generated  by $\{a_t: t \in K\}$. 
From our assumption it follows that
\[a_swa_s = u^{-1}v^{-1} \ \in \ A_J.\]
But the left hand side is an element of $A_{I \cup\{s\}}$. Now the 
associated partition must be non-crossing by Lemma \ref{lem:4}~{\bf (ii)}.
Hence $I,J$ and $\{s\}$ are pairwise disjoint
and $A_J \cap A_{I \cup\{s\}} = A_{\emptyset}$ by Theorem 1 of
\cite{appelschupp83}. We infer
\[ a_swa_s = \id.\]
\par
Repeating the argument leads to
\[ w \in A_{\{s\}} \mbox{ and } a_s^2 \in A_J.\]
The latter is in contradiction to $A_{\{s\}} \cap A_J = \{\id\}$. 
\par
The converse is true in the free group generated by $\{a_s : s\in S\}$,
hence in the Artin group.
\end{proof}
~\\
\begin{proofth}{\ref{artinclt}}
We compute as in the Coxeter case
\begin{eqnarray*}
\lefteqn{ \varphi_A \left( \left(\frac{\delta(a_1)+\delta(a_1^{-1})+ \ldots + %
\delta(a_N)+\delta(a_N^{-1})}{\sqrt{2N}}%
\right)^k\right)%
\ =}& &\\%
&=&\sum_{i_1, \ldots, i_k =1}^N \left(\frac{1}{\sqrt{2N}}\right)^k%
\sum_{\epsilon_1, \ldots, \epsilon_k \in \{+1,-1\}}%
\varphi_A \left(a_{i_1}^{\epsilon_1}\cdot \ldots \cdot a_{i_k}^{\epsilon_k}\right)\\
&=& \sum_{p=1}^k (2N)^{-\frac{k}{2}}%
\sum_{V\in \mathcal{V}_p^k}%
\sum_{w\in \mathcal{W}_{V}} \varphi_A ( \eval[w]),
\end{eqnarray*}
where $\mathcal{W}_{V}=%
\{(a_{i_1}^{\epsilon_1}, \ldots , a_{i_k}^{\epsilon_k}) : %
\mbox{  the word }(a_{i_1}^{\epsilon_1}, \ldots , a_{i_k}^{\epsilon_k})%
\mbox{ defines } V \}$.
For odd $k$ this vanishes, since in the Coxeter group
$\eval[\Phi(a_{i_1}^{\epsilon_1}),\ldots,\Phi(a_{i_k}^{\epsilon_k})] \neq
\id$.
\par
Further we assume $k=2n$ and we may omit in the summation partitions
which contain a one element set (Lemma \ref{lem:4}~{\bf (i)}).
The sum reducing to
\[ \sum_{p=1}^n (2N)^{-\frac{k}{2}}%
\sum_{V\in \mathcal{V}_p^k}%
\sum_{w\in \mathcal{W}_{V}} \varphi_A ( \eval[w]).
\]
Now, there are 
 $A_{N,p} = N (N-1) \cdot \ldots \cdot (N-p+1)$
choices of $a_{i_1},\ldots , a_{i_k}$ among 
$\{a_1,\ldots, a_N \}$ and $2^k$ choices of signs
$\epsilon_1,\ldots, \epsilon_k$ which define the same partition 
$V \in \mathcal{V}_p^{2n}$. For $p<n$ then
$(2N)^{-n} 2^k N (N-1) \cdot \ldots \cdot (N-p+1)$
tends to zero as $N$ tends to $\infty$.
With Lemma~\ref{lem:4} {\bf (ii)} we come to
\[
(2N)^{-n}\sum_{V\in \NC_2(2n)}%
\sum_{w\in \mathcal{W}_{V}} \varphi_A ( \eval[w]) + o(1).
\]
By Lemma~\ref{lem:5}, for any $V\in \NC_2(2n)$ we have 
$N (N-1)\cdot \ldots \cdot (N-n+1)$ choices of letters 
but only $2^n$ choices of signs, in building words $w\in \mathcal{W}_{V}$
which contribute a non-zero term in this sum. They then just contribute
$\varphi_A ( \eval[w])=1$.
In the limit we obtain
\[
\lim_{N\rightarrow \infty}(2N)^{-n}2^n\frac{N!}{n!} \card{\NC_2(k)} =%
\card{\NC_2(k)}
\]
\end{proofth}

{\bf author's address:}\\
Gero Fendler\\
Finstertal 16\\
D-69514 Laudenbach\\
Germany\\
\small
e--mail: gero.fendler@t-online.de

\begin{thebibliography}{1}

\bibitem{appelschupp83}
K.~Appel and P.~Schupp.
\newblock Artin groups and infinite {C}oxeter groups.
\newblock {\em Inventiones Math.}, 72:201--220, 1983.

\bibitem{bozspei91}
M.~Bo\.{z}ejko and R.~Speicher.
\newblock An example of generalized {B}rownian motion {I}.
\newblock {\em Communications in Math. Phys.}, 137:519--531, 1991.

\bibitem{bozspei92}
M.~Bo{\.{z}}ejko and R.~Speicher.
\newblock An example of generalized {B}rownian motion {II}.
\newblock In L.~Accardi, editor, {\em Quantum Probability and Related Topics
  {VII}.}, pages 67--77, Singapore, 1992. World Scientific.

\bibitem{Hum_90}
J.~E. Humphreys.
\newblock {\em {R}eflection {G}roups and {C}oxeter {G}roups.}
\newblock Cambridge University Press, Cambridge, 1990.

\bibitem{Kreweras72}
G.~Kreweras.
\newblock Sur les partitions non croisees d'un cycle.
\newblock {\em Discr. Math.}, 1:333--350, 1972.

\bibitem{speicher92}
R.~Speicher.
\newblock A non-commutative central limit theorem.
\newblock {\em Math. Z.}, 209:55--66, 1992.

\bibitem{szwarc_98}
R.~Szwarc.
\newblock Structure of geodesics in the {C}ayley graph of infinite {C}oxeter
  groups.
\newblock {\em Colloq. Math.}, 95:79--90, 2003.

\end{thebibliography}
\end{document}